# Convergence of Modified Picard-Mann Hybrid Iteration Process For Nearly Nonexpansive Mappings

Adrian Ghiura

*Department of Mathematics and Computer Science,*

*University Dunarea de Jos, Galati, Romania*

***Abstract*** *- In this paper, we prove the strong convergence theorems for nearly nonexpansive mappings, using the modified Picard-Mann hybrid iteration process in the context of uniformly convex Banach space.*

**2010 hematics Subject Classication:** 47H10, 47H09, 47J25.

## I. INTRODUCTION

In 1965, Browder [3] proved that every nonexpansive mapping in a uniformly convex Banach space has a fixed point. Goebel and Kirk [5] extended Browder's result to the class of asymptotically nonexpansive mappings. In 2005, Sahu [13] introduced the class of nearly Lipschitzian mappings which is an important generalization of the class of Lipschitzian mappings. Later, Agarwal et.al [1] proposed a new iteration process for the iterative approximation of fixed points of nearly asymptotically nonexpansive mappings. In 2013, Khan [8] introduced a new iteration process for nonexpansive mappings, which is called the Picard-Mann hybrid iteration process and showed that the new process converges faster than Picard and Mann iteration processes. Geethalakshmi and Hemavathy [4] proved strong convergence and stability results of the Picard-Mann hybrid iteration process for monotone nonexpansive mappings. Goyal [6] generalized a theorem of Xu, for weakly asymptotic contraction. Recently, Khan [9] proves the existence of fixed points of generalized nonexpansive mappings in CAT(0) spaces, and approximate them using Picard-Mann hybrid iteration processes. Akewe and Osilike [2] proved convergence and stability results for Picard-Mann hybrid iterative schemes for contractive-like operators in a real normed space.

Motived and inspired by this work, we introduce a modified Picard-Mann hybrid iteration process and prove some strong convergence theorems for nearly nonexpansive mappings in uniformly convex Banach space.

## II. PRELIMINARIES

Let $C$ be a nonempty subset of a real Banach space $X$ and $T: C \to C$ be a mapping with the fixed point set $F(T)$, i.e., $F(T) = \{p \in C : Tp = p\}$. Now, we recall some definitions and conclusions for our presentation.

**Definition 2.1.** ([5]) The function $\delta_X: [0,2] \to [0,1]$ is said to be the modulus of convexity of X if

$$\delta_X(\varepsilon) = \inf \left\{ 1 - \frac{\|x+y\|}{2} : \|x\| \leq 1, \|y\| \leq 1, \|x-y\| \geq \varepsilon \right\} \tag{2.1}$$

X is said to be uniformly convex if $\delta_X(0) = 0$ and $\delta_X(\varepsilon) > 0$ for all $\varepsilon \in (0,2]$.

**Definition 2.2.** A mapping T: $C \to C$ is said to be nonexpansive if

$$\|Tx - Ty\| \leq \|x - y\| \quad \text{for all } x, y \in C. \tag{2.2}$$

**Definition 2.3.** A mapping T: $C \to C$ is said to be asymptotically nonexpansive if there exists a sequence

$k_n \subset [1, \infty)$ with $\lim_{n \to \infty} k_n = 1$ such that

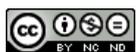 

$$\|T^n x - T^n y\| \leq k_n \|x - y\| \quad \text{for all} \ x, y \in C, \ n \geq 1 \tag{2.3}$$

**Definition 2.4.** A mapping $T: C \to C$ is said to be uniformly L-Lipschitzian if there exists a constant $L > 0$ such that

$$\|T^n x - T^n y\| \leq L \|x - y\| \quad \text{for all} \ x, y \in C, \ n \geq 1. \tag{2.4}$$

**Remark 2.1.** It easy to see that every nonexpansive mapping T is asymptotically nonexpansive with sequence $k_n = 1$ and every asymptotically nonexpansive mapping is uniformly L-Lipschitzian with $L = \sup_{n \in \mathbb{N}} k_n$.

**Definition 2.5.** ([13]) Let $\{a_n\}$ be a sequence in [0,1) with $\lim_{n \to \infty} a_n = 0$. A mapping $T: C \to C$ is said to be nearly nonexpansive with respect to $\{a_n\}$ if

$$\|T^n x - T^n y\| \leq \|x - y\| + a_n \quad \text{for all} \ x, y \in C, \ n \geq 1. \tag{2.5}$$

**Remark 2.2.** ([13]) If C is a bounded domain of an asymptotically nonexpansive mapping T, then T is nearly nonexpansive. In fact, we have

$$\|T^n x - T^n y\| \leq k_n \|x - y\|$$

$$\leq \|x - y\| + (k_n - 1)\|x - y\|$$

$$\leq \|x - y\| + (k_n - 1) \cdot diam(C), \quad \text{for all} \ x, y \in C, \ n \geq 1.$$

From Remark 2.1 and 2.2 we have the following implications:

nonexpansive $\Rightarrow$ asymptotically nonexpansive $\Rightarrow$ nearly nonexpansive

**Example 2.1.** Let $X = \mathbb{R}$, $C = [0,1]$ and $T: C \to C$ be a mapping defined by

$$Tx = \begin{cases} qx, & \text{if } x \in [0,1) \\ 0, & \text{if } x = 1. \end{cases}$$

where $q \in (0,1)$. It is clear that $T$ is discontinuous mapping. However, it is nearly nonexpansive mapping with respect to the sequence $a_n = q^n, a_n \to 0$. Indeed,

$$\|T^n x - T^n y\| \leq q^n \|x - y\| + q^n$$

$$\leq \|x - y\| + a_n \quad \text{for all} \ x, y \in C, \ n \geq 1.$$

**Lemma 2.1** ([11]). Let $\{a_n\}, \{b_n\}$ and $\{\delta_n\}$ be sequences of nonnegative real numbers satisfying the inequality

$$a_{n+1} \leq (1 + \delta_n) a_n + b_n, \quad \text{for all} \ n \in \mathbb{N}$$

If $\sum_{n=0}^{\infty} \delta_n < \infty$ and $\sum_{n=0}^{\infty} b_n < \infty$ then $\lim_{n \to \infty} a_n$ exists. In particular, if $\{a_n\}$ has a subsequence which converges strongly to zero, then $\lim_{n \to \infty} a_n = 0$.

**Lemma 2.2** ([14]). Let $X$ be a real uniformly convex Banach space and $0 < a \leq t_n \leq b < 1$, for all $n \in \mathbb{N}$. Let $\{x_n\}$ and $\{y_n\}$ be sequences in $X$ such that $\lim_{n \to \infty} \sup \|x_n\| \leq r$, $\lim_{n \to \infty} \sup \|y_n\| \leq r$, and $\lim_{n \to \infty} \|(1 - t_n)x_n + t_n y_n\| = r$, hold for some $r \geq 0$. Then $\lim_{n \to \infty} \|x_n - y_n\| = 0$.

We will now consider some well-known iteration schemes. Let $C$ be a nonempty convex subset of normed space $X$ and $T: C \to C$ a self-map.

(a) The Picard iteration process is defined by

$$x_{n+1} = Tx_n$$

for all $n \geq 0$, (see for more information [12]).

(b) The Mann iteration process (see, for example [10]) is defined by

$$x_{n+1} = (1 - \alpha_n)x_n + Tx_n$$

for all $n \geq 0$ and $\{\alpha_n\}_{n \geq 0}$ is a real sequence in $[0, 1]$ which satisfies the conditions: $0 \leq \alpha_n < 1$ and $\sum_{n=0}^{\infty} \alpha_n = \infty$.

(c) The Ishikawa iteration process (see, for example [7]) is defined by

$$x_{n+1} = (1 - \alpha_n)x_n + \alpha_n T y_n$$
$$y_n = (1 - \beta_n)x_n + \beta_n T x_n$$

where $\{\alpha_n\}_{n \geq 0}$ and $\{\beta_n\}_{n \geq 0}$ be real sequences in $[0, 1]$ for all $n \geq 0$.

(d) The modified Mann iteration process (see, for example [14]) is defined by

$$x_{n+1} = (1 - \alpha_n)x_n + T^n x_n$$

where $\{\alpha_n\}_{n \geq 0}$ is a real sequence in $[0, 1]$ which satisfies condition $0 < a \leq \alpha_n \leq b < 1$ for all $n \geq 0$.

(e) The Picard-Mann hybrid iteration process (see, for example [8]) is defined by

$$x_{n+1} = T y_n$$
$$y_n = (1 - \alpha_n)x_n + \alpha_n T x_n$$

for all $n \geq 0$ and $\{\alpha_n\}_{n \geq 0}$ is a real sequence in $[0, 1]$.

### III. MAIN RESULTS

We introduce a modified Picard-Mann hybrid iteration process by

$$\begin{aligned} x_{n+1} &= T^n y_n \\ y_n &= (1 - \alpha_n)x_n + \alpha_n T^n x_n, \end{aligned} \quad (3.1)$$

for all $n \geq 0$ and $\{\alpha_n\}_{n \geq 0}$ is a real sequence in $[0, 1]$.

The purpose of this section is to prove some strong convergence theorems with respect to the iteration scheme (3.1) for nearly nonexpansive mapping in real uniformly convex Banach spaces.

**Theorem 3.1.** Let $C$ be a nonempty compact convex subset of a real uniformly convex Banach space X. Let $T: C \to C$ be a uniformly L-Lipschitzian, nearly nonexpansive mapping with respect to $\{a_n\}$ such that $\sum_{n=0}^{\infty} a_n < \infty$. Let $\{x_n\}$ be the modified Picard-Mann iteration defined by (3.1) where $\{\alpha_n\}$ is sequence in $(0,1)$. Then the following hold:

*(i)* $\lim_{n\to\infty}\|x_n - p\|$ exists for all $p \in F(T)$;

*(ii)* $\lim_{n\to\infty}\|x_n - T^n x_n\| = 0$;

*(iii)* $\{x_n\}$ converges strongly to a fixed point of T.

*Proof. (i)* By Schauder's fixed point theorem, we obtain that $F(T) \neq \emptyset$. Let $p \in F(T)$, by (3.1) we have

$$\|x_{n+1} - p\| = \|T^n x_n - p\| \leq \|y_n - p\| + a_n$$

$$\|y_n - p\| = (1 - \alpha_n)\|x_n - p\| + \alpha_n\|T^n x_n - p\|$$
$$\leq (1 - \alpha_n)\|x_n - p\| + \alpha_n(\|x_n - p\| + a_n)$$
$$= \|x_n - p\| + (1 + \alpha_n)a_n \qquad (3.2)$$

So, we have

$$\|x_{n+1} - p\| \leq \|x_n - p\| + (1 + \alpha_n)a_n \qquad (3.3)$$

It follows from $\sum_{n=0}^{\infty} a_n < \infty$ and Lemma 2.1 that $\lim_{n\to\infty}\|x_n - p\|$ exists for $p \in F(T)$.

*(ii)* We set $\lim_{n\to\infty}\|x_n - p\| = c$, from (3.2) we have

$$\lim_{n\to\infty} \sup\|y_n - p\| \leq c \qquad (3.4)$$

Also,

$$\|T^n y_n - p\| \leq \|y_n - p\| + a_n$$

So, we have

$$\lim_{n\to\infty} \sup\|T^n y_n - p\| \leq c \qquad (3.5)$$

Similarly

$$\|T^n x_n - p\| \leq \|x_n - p\| + a_n$$

and we have

$$\lim_{n\to\infty} \sup\|T^n x_n - p\| \leq c \qquad (3.6)$$

Now

$$\|x_{n+1} - p\| \leq \|y_n - p\| + a_n$$

Taking limit infimum, we have

$$\lim_{n\to\infty} \inf \|x_{n+1} - p\| \leq \lim_{n\to\infty} \inf \|y_n - p\|$$

$$c \leq \lim_{n\to\infty} \inf \|y_n - p\| \leq c \qquad (3.7)$$

From (3.4) and (3.7) we have

$$c = \lim_{n\to\infty} \|y_n - p\| \tag{3.8}$$

i.e.,

$$c = \lim_{n\to\infty} \|y_n - p\| = \lim_{n\to\infty} \|(1 - \alpha_n)(x_n - p) + \alpha_n(T^n x_n - p)\|$$

Therefore by Lemma 2.2, we obtain

$$\lim_{n\to\infty} \|x_n - T^n x_n\| = 0 . \tag{3.9}$$

*(iii)* Finaly, since $T$ is uniformly $L$-Lipschitzian mapping, we have

$$\|x_n - Tx_n\| \leq \|x_n - x_{n+1}\| + \|x_{n+1} - T^{n+1}x_{n+1}\| + \|T^{n+1}x_{n+1} - T^{n+1}x_n\| + \|T^{n+1}x_n - Tx_n\|$$

$$\leq \|x_n - x_{n+1}\| + \|x_{n+1} - T^{n+1}x_{n+1}\| + L\|x_{n+1} - x_n\| + L\|x_n - T^n x_n\|$$

Since $T$ is uniformly continuous, it follows from (3.9) that

$$\lim_{n\to\infty} \|x_n - Tx_n\| = 0 . \tag{3.10}$$

By the compactness of $C$, there is a subsequence $\{x_{n_k}\}$ of $\{x_n\}$ such that

$$\lim_{n\to\infty} x_{n_k} = p^* \tag{3.11}$$

Since $T$ is continuous, it follows from (3.10) that $p^* \in F(T)$. Since $\lim_{n\to\infty} \|x_n - p\|$ exists for all $p \in F(T)$, and we conclude from (3.11) that $\lim_{n\to\infty} x_n = p^* \in F(T)$.

**Theorem 3.2.** Assume that all the conditions of Theorem 3.1 are satisfied. Then the sequence $\{x_n\}$ generated by (3.1) converges strongly to a fixed point of $T$ if and only if

$$\lim_{n\to\infty} \inf d(x_n, F(T)) = 0,$$

where $d(x, F(T)) = \inf \{d(x,p) : p \in F(T)\}$.

*Proof.* Necessity is obvious. Conversely, suppose that $\lim_{n\to\infty} \inf d(x_n, F(T)) = 0$. From Theorem 3.1 we know that $\lim_{n\to\infty} \|x_n - p\|$ exists for all $p \in F(T)$, so $\lim_{n\to\infty} d(x_n, F(T))$ exists for all $p \in F(T)$. Thus by hypothesis

$$\lim_{n\to\infty} d(x_n, F(T)) = 0. \tag{3.12}$$

Now we show that $\{x_n\}$ is a Cauchy sequence in $C$. Indeed, from (3.3), we have

$$\|x_{n+1} - p\| \leq \|x_n - p\| + (1 + \alpha_n)a_n$$

Now, we set $b_n := (1 + \alpha_n)a_n$. For any $m, n \in \mathbb{N}$, $m > n \geq 1$, we have

$$\|x_m - p\| \leq \|x_{m-1} - p\| + b_{m-1}$$
$$\leq \|x_{m-2} - p\| + b_{m-1} + b_{m-2}$$
$$\vdots$$
$$\leq \|x_n - p\| + \sum_{i=0}^{m-1} b_i$$
$$\leq \|x_n - p\| + \sum_{i=0}^{\infty} b_i$$

Since $\lim_{n\to\infty} d(x_n, F(T)) = 0$ and $\sum_{i=0}^{\infty} b_i < \infty$ for any $\varepsilon > 0$ there exists a positive integer $n_0$ such that

$$d(x_n, F(T)) < \frac{\varepsilon}{4}, \quad \sum_{i=n_0}^{\infty} b_i < \frac{\varepsilon}{4}.$$

Therefore, there exists $\bar{p} \in F(T)$ such that

$$\|x_{n_0} - \bar{p}\| < \frac{\varepsilon}{4}, \quad \sum_{i=n_0}^{\infty} b_i < \frac{\varepsilon}{4}.$$

Thus, for all $m, n \geq n_0$ we get from the above inequality that

$$\|x_m - x_n\| \leq \|x_m - \bar{p}\| + \|x_n - \bar{p}\|$$

$$\leq \|x_{n_0} - \bar{p}\| + \sum_{i=n_0}^{\infty} b_i + \|x_{n_0} - \bar{p}\| + \sum_{i=n_0}^{\infty} b_i$$

$$= 2\left(\|x_{n_0} - \bar{p}\| + \sum_{i=n_0}^{\infty} b_i\right)$$

$$< 2\left(\frac{\varepsilon}{4} + \frac{\varepsilon}{4}\right) = \varepsilon.$$

Thus, it follows that $\{x_n\}$ is a Cauchy sequence. Since $C$ is a closed subset of Banach space $X$, the sequence $\{x_n\}$ converges strongly to some $p^* \in C$. Since $F(T)$ is a closed subset of $C$ and $\lim_{n\to\infty} d(x_n, F(T)) = 0$ we have $p^* \in F(T)$. Thus, the sequence $\{x_n\}$ converges strongly to a fixed point of $T$. This completes the proof.

Senter and Dotson [15] introduced the notion of mapping satisfying Condition (I) which is defined as follows:

**Definition 3.1.** A mapping $T: C \to C$ is said to satisfy Condition (I), if there exists a non-decreasing function

$\varphi : [0, \infty) \to [0, \infty)$ with $\varphi(0) = 0$ and $\varphi(t) > 0$, for all $t > 0$ such that

$$d(x, Tx) \geq \varphi(d(x, F(T)))$$

for all $x \in C$.

**Theorem 3.3.** Assume that all the conditions of Theorem 3.1 are satisfied and let $T$ be a mapping satisfying Condition (I). Then the sequence $\{x_n\}$ generated by (3.1) converges strongly to a fixed point of $T$.

*Proof.* We proved in Theorem 3.2, that $\lim_{n\to\infty} d(x_n, F(T))$ exists. From Theorem 3.1 we have $\lim_{n\to\infty} d(x_n, Tx_n) = 0$.

It follows from Condition (I) that

$$\lim_{n\to\infty} \varphi\left(d(x_n, F(T))\right) \leq \lim_{n\to\infty} d(x_n, Tx_n) = 0$$

i.e,
$$\lim_{n\to\infty} \varphi\left(d(x_n, F(T))\right) = 0.$$

Since $\varphi : [0, \infty) \to [0, \infty)$ is a non-decreasing function with $\varphi(0) = 0$ and $\varphi(t) > 0$, for all $t > 0$, we obtain

$$\lim_{n\to\infty} d(x_n, F(T)) = 0.$$

Consequently, $\{x_n\}$ converges strongly to a fixed point of $T$.

## REFERENCES

[1] R. P. Agarwal, D. O'Regan and D. R. Sahu, Iterative construction of fixed points of nearly asymptotically nonexpansive mappings, J. Nonlinear Convex Anal., 8 (1) (2007), 61-79.

[2] H. Akewe, G. A. Okeke, Convergence and stability theorems for the Picard-Mann hybrid iterative scheme for a general class of contractive-like operators, Fixed Point Theory Appl., 66 (2015), 1-8.

[3] F. E. Browder, Nonexpansive nonlinear operators in Banach space, Proc. Nat. Acad. Sci. U.S.A., 54 (1965), 1041-1044.

[4] B. Geethalakshmi, R. Hemavathy, Picard-Mann Hybrid Iteration Processes for Monotone Non Expansive Mappings, International Journal of Mathematics Trends and Technology, 55 (5) (2018), 367-370.

[5] K. Goebel and W. A. Kirk, A fixed point theorem of asymptotically nonexpansive mappings, Proc. Amer. Math. Soc., 35 (1972), 171-174.

[6] S. Goyal, A Generalization of a fixed point theorem of Hong-Kun Xu, International Journal of Mathematics Trends and Technology, 38 (1) (2016), 23-26.

[7] S. Ishikawa, Fixed points by a new iteration method, Proc.Am. Math. Soc., 44 (1974), 147-150.

[8] S. H. Khan, A Picard-Mann hybrid iterative process, Fixed Point Theory Appl., 69 (2013), 1-10.

[9] Ritika and S. H. Khan, Convergence of Picard-Mann Hybrid Iterative Process for Generalized Nonexpansive Mappings in CAT(0) Spaces, Filomat, 31 (11) (2017), 3531-3538.

[10] W. R. Mann, Mean value methods in iteration, Proc. Am. Math. Soc., 4 (1953), 506-510.

[11] M. O. Osilike and S. C. Aniagbosor, Weak and strong convergence theorems for fixed points of asymptotically nonexpansive mappings, Math. Comput. Modelling, 32 (2000), 1181-1191.

[12] E. Picard, Memoire sur la theorie des equations aux derivee partielles et la metode des approximations successives, J. Math. Pures. Appl., 6 (1890), 145-210.

[13] D. R. Sahu, Fixed points of demicontinuous nearly Lipschitzian mappings in Banach spaces, Comment. Math. Univ. Carolin, 46 (2005), 653-666.

[14] J. Schu, Iterative construction of fixed points of asymptotically nonexpansive mappings, J. Math. Anal. Appl., 158 (1991), 407-412.

[15] H. F. Senter and W. G. Dotson, Approximating fixed points of nonexpansive mappings, Proc. Am. Math. Soc., 44 (2) (1974), 375-380.